\documentclass{article}

\usepackage{amsfonts}
\usepackage[utf8]{inputenc}
\begin{document}

\title{The Geometry of $H_4$ Polytopes}

\author{Tomme Denney, Da'Shay Hooker, De'Janeke Johnson,\\ {Tianna Robinson, Majid Butler and Sandernisha Claiborne}{
\renewcommand{\thefootnote}{\fnsymbol{footnote}}
\footnote{The authors were supported by the Student Research 
Fellowship program at McDonogh 35 High School in New Orleans. The authors thank Professor Cathy Kriloff of Idaho State University for suggesting they use the techniques they developed in [1] to study $H_3$ (and $H_4$) symmetry, they thank their advisor Rich Margolin for several helpful conversations, and they thank the referee for a number of suggestions that improved the text and corrected the terminology of this paper.}
\setcounter{footnote}{0}
}
\\{ }
\\{ }
 \\
{ }}

\date{We describe the geometry of an arrangement of 24-cells inscribed in the 600-cell. In $\S$7 we apply our results to the even unimodular lattice $E_8$ and show how the 600-cell transforms $E_8$/2$E_8$, an 8-space over the field $\bf{F}$$_2$, into a 4-space over $\bf{F}$$_4$ whose points, lines and planes are labeled by the geometric objects of the 600-cell.\\
\leavevmode\\
May 25, 2019}
\maketitle

\section{Introduction}

The 600-cell occupies a difficult niche in Euclidean geometry. 
Because it exists in four dimensions the Greeks never considered it, and it will never achieve the fame of the five Platonic solids. 
Because it’s the last of its kind (dimensions higher than four admit only simplexes, orthoplexes and n-cubes ([14], Chapter X)), no higher-dimensional polytope analogues of the 600-cell exist, so it can never be part of a systematic treatment of polytopes and must always be treated as an outlier. 
Even its most natural coordinates in 4-space are a disadvantage: they lie in $\mathbb{Z}$[$\frac{1+\sqrt{5}}{2}$] and not in $\mathbb{Z}$, so the regular $H_4$ polytopes (the 600-cell and its dual, the 120-cell) seem less natural than polytopes whose coordinates are all integers.

And yet these drawbacks account for the 600-cell’s fascination. 
The 600-cell is scarcely known because it exists in four dimensions and not three; but a well-known theorem of Hurwitz states that composition algebras exist only in dimensions 1, 2, 4 and 8—--these are the real numbers, the complex numbers, the quaternions and the octonions ([9], $\S$7.6)—--which raises the possibility that the 600-cell has some connection to the quaternions. 
(In fact, it does: see Fact 8 and $\S$2 below).
And because the four coordinates for $H_4$ polytopes lie in $\mathbb{Q}$($\sqrt{5}$), the variety of “norm reduction” maps from this extension field to $\mathbb{Q}$ results in numerous different ways to view these polytopes ($\S$6, Examples 1-3). 
Still, despite these varied approaches to $H_4$ objects, the geometry of the 600-cell (and other $H_4$ polytopes) has remained all too mysterious.

Even the experts have sometimes been confused about the structure of $H_4$ polytopes.
In 1907 the Dutch mathematician P.H. Schoute published, without proof, the following claim: “Die 120 Eckpunkte eines Z$^{(1)}_{600}$ bilden auf zehn verschiedene Arten die Eckpunkte von f\"{u}nf Z$^{\frac{1}{2}(1+\sqrt{5})}_{24}$.” ([13], p. 231) (“The 120 vertices of a 600-cell form in ten different ways the vertices of five 24-cells.”) 
Evidently this claim was viewed skeptically: in 1933 H.S.M. Coxeter published his opinion that “But surely \textit{zehn} should be \textit{f\"{u}nf} in the phrase \textit{auf zehn verschiedene Arten},” ([4], p. 337) i.e., that surely there were only five and not ten such partitions. 

But in fact there are ten, and not five, such partitions of the vertices of the 600-cell.
Coxeter realized his mistake and corrected it: “Thus Schoute was right when he said the 120 vertices of $\{$3, 3, 5$\}$ belong to five $\{$3, 4, 3$\}$’s in ten different ways. 
The disparaging remark in the second footnote to Coxeter [4], p. 337, should be deleted.” ([5], p. 270)
And although it was later fixed, Coxeter’s mistake highlights a difficulty that has plagued $H_4$ polytopes (and their symmetry group) since their study began roughly 150 years ago: $H_4$ remains more complicated than it should.

So the goal of our paper is to make the geometry of the 600-cell a little more accessible.
We've tried to achieve this in four ways: (i) by collecting a number of results (some old and some new) about the 600-cell and listing them as Facts later in our Introduction; (ii) by proving Schoute's century-old statement that there are exactly ten ways to partition the 120 vertices of a 600-cell into five disjoint 24-cells\footnote{Throughout this paper we identify a polytope with its vertex set.} ([13], p. 231 and $\S$3); (iii) by using a refinement of these ten partitions to label the vertices of the 600-cell and the 120-cell so that $H_4$ symmetry becomes more transparent ($\S$2); and (iv) by providing new ways to look at the 600-cell, in the hope of finding even more avenues towards understanding $H_4$ geometry ($\S$$\S$6 and 7 and also [10].)

\bigskip
\medskip

We begin with our ten Facts:

\medskip

\textbf{Fact 1:} The 600-cell is a regular $H_4$ polytope that consists of 120 vertices, 720 edges, 1200 triangular faces and 600 tetrahedral cells. ([12], p. 386; [5], p. 153)

\medskip

\textbf{Fact 2:} The 600-cell contains exactly 25 24-cells, 75 16-cells and 75 8-cells, with 
each 16-cell and each 8-cell lying in just one 24-cell. ($\S$3)

\medskip

\textbf{Fact 3:} The rotation group of the 600-cell is the central product 2$A_5$ $\circ$ 2$A_5$ $\cong$ 2($A_5$ $\times$ $A_5$) ([11], p. 17) so the full symmetry group (including reflections) of a 600-cell $H$ is Aut($H$) $\cong$ 2($A_5$ $\times$ $A_5$)2.

\medskip

\textbf{Fact 4:} The subgroup of Aut($H$) that fixes a vertex is isomorphic to 2 $\times$ $A_5$ and the subgroup that stabilizes a 24-cell inscribed into a 600-cell is isomorphic to 2($A_4$ $\times$ $A_4$)2. ($\S$7)

\medskip

\textbf{Fact 5:} The 25 24-cells can be placed in a 5 $\times$ 5 array, so that each row and each column of the array partition the 120 vertices of the 600-cell into five disjoint 24-cells. 
The rows and columns of the array are the only ten such partitions of the 600-cell. ([13], p. 231, stated without proof, and $\S$$\S$2, 3)

\medskip

\textbf{Fact 6:} The quotient group Aut($H$)/$\{$$\pm$1$\}$ $\cong$ ($A_5$ $\times$ $A_5$)2 is isomorphic to a subgroup of the symmetric group $S_{10}$, where the first $A_5$ permutes five symbols 1, 2, 3, 4 and 5 and fixes five symbols 6, 7, 8, 9 and X; the second $A_5$  permutes 6, 7, 8, 9 and X and fixes 1, 2, 3, 4 and 5; and elements outside the subquotient $A_5$ $\times$ $A_5$ interchange the pentads $\{$1, 2, 3, 4, 5$\}$ and $\{$6, 7, 8, 9, X$\}$. ($\S\S$2, 3)

\medskip

\textbf{Fact 7:} The 25 24-cells of a 600-cell may be labeled by duads ($i$ $j$), where 1 $\leq$ $i$ $\leq$ 5 and 6 $\leq$ $j$ $\leq$ X, and the $\pm$60 vertices\footnote{By "$\pm$60 vertices" we mean the 60 pairs $\{$$v, -v$$\}$ of opposite vertices.} of the 600-cell may be labeled by five duads (1 $j_1$)(2 $j_2$)(3 $j_3$)(4 $j_4$)(5 $j_5$) where ($j_1$, $j_2$, $j_3$, $j_4$, $j_5$) is an even permutation of (6, 7, 8, 9, X). 
The symmetry group Aut($H$) permutes the 24-cells and the vertices of the 600-cell as ($A_5$ $\times$ $A_5$)2 $\cong$ Aut($H$)/$\{$$\pm$1$\}$ permutes their labels. ($\S$2)

\medskip

\textbf{Fact 8:} The 120 vertices of a 600-cell may be viewed as the 120 \textit{icosians}, the units of a quaternion algebra whose coefficients lie in $\mathbb{Q}$($\sqrt{5}$). ([6], p. 74)

\medskip

\textbf{Fact 9:} A 600-cell $H$ can be "embedded" in the $E_8$ lattice by means of a "norm reduction" map, so that the 120 vertices of $H$ become 120 of the 240 root vectors of $E_8$. 
The remaining 120 root vectors are the vertices of the embedded 600-cell $\varphi$$H$ (a scale copy of $H$), where $\varphi$ is the golden ratio $\frac{1}{2}$(1+$\sqrt{5}$). ([15], p. 59)

\medskip

\textbf{Fact 10:} The "embeddings" of the 600-cells $H$ and $\varphi$$H$ in $E_8$ turn $E_8$/2$E_8$ (an 8-space over $\bf{F}$$_2$) into a 4-space over $\bf{F}$$_4$ whose 85 points are labeled by the $\pm$60 vertices and the 25 24-cells of $H$. ($\S$7)

\bigskip

In addition, a few other aspects of $H_4$ geometry (of lesser importance) are italicized in $\S$$\S$2-4.

\section{The Labeling of $H_4$}

In this section (much of which follows $\S$4.27 of [6]) we exploit the algebraic properties of the 600-cell that are evident when we view its 120 vertices as the 120 icosians—that is, when we consider the vertices\hfill

\hspace{2.3cm}($\pm$1, 0, 0, 0)$^S$, \hspace{1.3cm} (4 pairs of vertices)\hfill

\hspace{2.1cm}$\frac{1}{2}$($\pm$1, $\pm$1, $\pm$1, $\pm$1), \hspace{0.75cm}    (8)\hfill

\hspace{2.1cm}$\frac{1}{2}$(0, $\pm$1, $\pm$$\varphi$, $\pm$$\varphi^{-1}$)$^A$ \hspace{0.35cm} (48) \hfill\break
as a set of 120 quaternions closed under multiplication, where the four coordinates are scalars for (1, $i$, $j$, $k$) and where $\varphi$ = $\frac{1}{2}$($1+$$\sqrt{5}$), ${\varphi}^{-1}$ = $\frac{1}{2}$($-1+$$\sqrt{5}$), $S$ $\equiv$ all permutations of the coordinates and $A$ $\equiv$ all even permutations are permitted. ([6], p. 74; [11], p. 16)

We will let 2$A_5$ denote both the set of 120 vertices and also the symmetry subgroup these vertices generate, and we will sometimes write 2$A_{5L}$ or 2$A_{5R}$ to emphasize that 2$A_5$ is acting by left or right icosian multiplication, respectively.
Similarly 2$A_4$ will denote both the sub-polytope of 24 vertices (a 24-cell with 8 vertices ($\pm$1, 0, 0, 0)$^S$ and 16 vertices $\frac{1}{2}$($\pm$1, $\pm$1, $\pm$1, $\pm$1)) and also the subgroup these vertices generate, with 2$A_{4L}$ or 2$A_{4R}$ specifically denoting left or right multiplication.
It will be clear from the context which meaning is intended.

We begin with the 24-cell 2$A_4$. If $g$ is an icosian of order 5 then $\{$$g^i$$\}$ (0 $\leq$ $i$ $\leq$ 4) is a left transversal of 2$A_4$ in 2$A_5$ and  $\{$$g^{-j}$$\}$ (0 $\leq$ $j$ $\leq$ 4) is a right transversal. 
Since 2$A_4$ is stable under both 2$A_{4L}$ and 2$A_{4R}$ and since 2$A_4$ is self-normalizing in 2$A_5$, we obtain 25 distinct 24-cells $\{$$g^i$2$A_4$$g^{-j}$$\}$$_{0 \leq i, j \leq 4}$ under left and right multiplication by 2$A_5$; that is, the rotation subgroup 2($A_5$$\times$$A_5$) of Aut($H$) yields 25 different 24-cells. 
In the next section we prove that these 25 are the only 24-cells within a 600-cell.

These 25 24-cells may be arranged in a 5$\times$5 array (\textbf{Fact 5}) by letting the ($i$, $j$)$^{th}$ entry be the 24-cell $g^i$2$A_4$$g^{-j}$, so that $A_{5L}$ permutes the rows of the array and $A_{5R}$ permutes the columns. 
Since the cosets of 2$A_4$ in 2$A_5$ contain each element once, the five 24-cells $\{$$g^i$2$A_4$$\}$ (the first column of the array) form one of Schoute’s partitions of the 600-cell; since right multiplication by any right coset representative $g^{-j}$ yields a symmetry of the 600-cell, each column of the array is also a partition. 
But the same arguments apply to the top row $\{$2$A_4$$g^{-j}$$\}$ of the array and to each of the five rows in general. 
Thus we obtain ten partitions of the 600-cell into disjoint 24-cells---the five rows and five columns of the array---and in the next section we prove that these ten are the only ones possible.

In the meantime we'll use the array and its ten partitions to label the 60 vertex pairs and the 25 24-cells of the 600-cell, creating a simple way to calculate in the quotient group Aut($H$)/$\{$$\pm$1$\}$ $\cong$ ($A_5$$\times$$A_5$)2. 
First, we label the five rows of the array by the symbols 1, 2, 3, 4 and 5 and the five columns by 6, 7, 8, 9 and X, so that each of Schoute’s partitions corresponds to one of the ten symbols. 
Next, we label each of the 25 24-cells in the array by a duad ($i$ $j$), where 1 $\leq$ $i$ $\leq$ 5 and 6 $\leq$ $j$ $\leq$ X, according to its position (\textbf{Fact 7}). 

Furthermore, since a 24-cell contains 12 vertex pairs $\{$$\pm$$v$$\}$, each of the 60 vertex pairs of the 600-cell must lie in five of the 25 24-cells, so we can name a vertex (pair) by the five different 24-cells that contain it.
Since all $\pm$60 vertices of the 600-cell appear exactly once in each row and once in each column, each vertex (pair) of the 600-cell is labeled by five duads (1 $j_1$)(2 $j_2$)(3 $j_3$)(4 $j_4$)(5 $j_5$) where ($j_1$, $j_2$, $j_3$, $j_4$, $j_5$) is some permutation of (6, 7, 8, 9, X).
The unit icosian 1$_I$ = (1, 0, 0, 0) is labeled (16)(27)(38)(49)(5X), as the 24-cells along the main diagonal of the array are the five 24-cells $g^i$2$A_4$$g^{-i}$ that contain $g^i$$1_I$$g^{-i}$ = 1$_I$. 

If a vertex $v$ has label (1 $j_1$)(2 $j_2$)(3 $j_3$)(4 $j_4$)(5 $j_5$), then $v_R$ (right multiplication by the icosian $v$) must send partition 6 to partition $j_1$, and sends 7 to $j_2$, etc., since $v_R$ sends the 24-cell (16) that contains the unit icosian 1$_I$ to the 24-cell in row 1 that contains $v$, namely, (1 $j_1$). 
Thus $v_R$ induces an even permutation of $\{$6, 7, 8, 9, X$\}$.  
Alternatively, the 24-cells that contain $\{$$\pm$$v$$\}$ can instead demonstrate how the rows are permuted by $v_L$: if the same five duads (1 $j_1$)(2 $j_2$)(3 $j_3$)(4 $j_4$)(5 $j_5$) are now written (in a different order, perhaps) as ($i_1$ 6)($i_2$ 7)($i_3$ 8)($i_4$ 9)($i_5$ X), then evidently $v_L$ sends partition 1 to partition $i_1$, sends 2 to $i_2$, etc.
Thus $v_L$ induces an even permutation of $\{$1, 2, 3, 4, 5$\}$, and so \emph{the permutations induced by $v_R$ and $v_L$ are indicated by the label of $v$.}

\smallskip

The 60 reflections that generate Aut($H$) are also easy to describe using the ten symbols. 
In the next section we will show that elements in the odd half of Aut($H$)—--the elements outside the rotation subgroup 2($A_5$$\times$$A_5$) of index 2---interchange the rows and columns of our 5$\times$5 array. 
In other words, each of these elements (including the 60 reflections) interchanges the symbols $\{$1, 2, 3, 4, 5$\}$ with $\{$6, 7, 8, 9, X$\}$. 
But for any vertex $v$ of the 600-cell, the reflection $r_v$ (that uses $v$$^\perp$ as its hyperplane of reflection) stabilizes each 24-cell that contains $v$, since 2$A_4$ is itself a closed reflection group; and because the five 24-cells {(1 $j_1$), (2 $j_2$), (3 $j_3$), (4 $j_4$), (5 $j_5$)} maintain their positions in the 5$\times$5 array, \emph{the reflection $r_v$ permutes Schoute’s ten partitions as (1 $j_1$)(2 $j_2$)(3 $j_3$)(4 $j_4$)(5 $j_5$)}. 
Hence the label of a vertex pair $\{$$\pm$$v$$\}$ displays the permutation induced by $r_v$, and the entire quotient group Aut($H$)/$\{$$\pm$1$\}$ $\cong$ ($A_5$$\times$$A_5$)2 obtains a straightforward description as a transitive but imprimitive subgroup of $S_{10}$ that stabilizes the two pentads $\{$1, 2, 3, 4, 5$\}$ and $\{$6, 7, 8, 9, X$\}$ (\textbf{Fact 6}).

\medskip

A similar process labels the 300 pairs of vertices of the 120-cell. 
We start with the tetrahedral cell of the 600-cell that has vertices 

\medskip

\hspace{2.5cm}$\frac{1}{2}$($\varphi$, 1, $\varphi$$^{-1}$, 0), \hspace{1.3cm} \hfill

\hspace{2.5cm}$\frac{1}{2}$($\varphi$, 1, -$\varphi$$^{-1}$, 0), \hspace{0.75cm}    \hfill

\hspace{2.5cm}$\frac{1}{2}$(1, $\varphi$, 0, $\varphi$$^{-1}$) and \hspace{0.75cm}    \hfill

\hspace{2.5cm}$\frac{1}{2}$(1, $\varphi$, 0, -$\varphi$$^{-1}$).\hspace{0.35cm} \hfill

\medskip
The center $\frac{1}{4}$($\varphi$+1, $\varphi$+1, 0, 0) of this tetrahedron is then a vertex of the dual of the 600-cell---that is, a vertex of the 120-cell. 
If we rescale we obtain the standard coordinates (2, 2, 0, 0) of a vertex $v_1$ of the 120-cell, and if we rotate $v_1$ by elements of 2$A_4$ (acting on the left and/or the right) we obtain a 24-cell $C$ = 2$A_4$$v_1$2$A_4$ = $\{$($\pm$2, $\pm$2, 0, 0)$^S$$\}$ as a sub-polytope of the 120-cell. 

As before, we obtain 25 distinct 24-cells $\{$$g^i$$C$$g^{-j}$$\}$ permuted by Aut($H$) in a 5$\times$5 array; since 2($A_5$$\times$$A_5$) is transitive on the 600 vertices of the 120-cell these 25 24-cells are mutually disjoint. 
Hence each 24-cell (and each of its 12 pairs of vertices) can initially be labeled by a duad ($i$ $j$), where 1 $\leq$ $i$ $\leq$ 5 and 6 $\leq$ $j$ $\leq$ X.
If we further label each vertex (pair) by appending the label of each of its four neighbors--—the four other vertices with which it shares an edge--—then \emph{each of the $\pm$300 vertices of the 120-cell is uniquely labeled by 1+4 duads}. 
For example, (38)$\big{|}$(16)(27)(4X)(59) denotes the vertex in the 24-cell (38) whose four neighbors lie in the 24-cells $\{$(16), (27), (4X), (59)$\}$. 
Naturally the labels for the vertices of the 120-cell are "duals" of the labels for the vertices of the 600-cell, in the sense that our example's duads $\{$(38), (16), (27), (4X), (59)$\}$ are transformed to $\{$(16), (27), (38), (49), (5X)$\}$ (the duads for the icosian unit 1$_I$ of the 600-cell) by an odd permutation of $\{$6, ..., X$\}$, whereas labels for vertices of the 600-cell are conjugate under even permutations.

\smallskip

Finally, since the first column of the 120-cell's array is the orbit (under icosian multiplication 2$A_{5L}$) of a single vertex of the 120-cell, the 120 vertices in the first column obtain the geometry of the icosians, i.e., they form a 600-cell. 
The same is true for each of the four other columns, which shows that \emph{the 600 vertices of the 120-cell may be partitioned into five disjoint 600-cells}.
Likewise, by using 2$A_{5R}$ we find that the five rows partition the 120-cell into disjoint 600-cells as well.

\section{An Overdue Proof of Schoute's Result}

In this section we prove Schoute's stated result that there are exactly ten ways to partition the 120 vertices of a 600-cell into five disjoint 24-cells ([13], p. 231).
To keep things simple, we'll use the standard coordinates for the 600-cell (i.e., twice the icosian coordinates of $\S$2) where the inner products of vertices are the natural inner product divided by 2.

\medskip

We begin by proving \textbf{Fact 2}.\hfill\break

$\textbf{Theorem}$: The 600-cell contains exactly 25 24-cells, 75 16-cells and 75 8-cells, with each 16-cell and each 8-cell lying in just one 24-cell. 

\smallskip Proof: To begin with, 2($A_5$$\times$$A_5$)2 is transitive on the 120 vertices of the 600-cell, since given any two vertices $v$ and $w$, icosian multiplication by ($v^{-1}$$w$)$_R$ $\in$ 2$A_{5R}$ sends $v$ to $w$.
Furthermore, given a vertex $v$ its stabilizer Stab($v$) is transitive on the set of all vertices that have a given inner product with $v$; for example, if we let $v$ = $1_I$ = (2, 0, 0, 0) then the monomial subgroup $2^3$:$A_3$ $\leq$ Stab($v$) (which leaves the first coordinate untouched) is itself transitive on the set $\mathcal{S_\varphi}$ of 12 vertices that have inner product 
$\varphi$ with $v$, and $2^3$:$A_3$ is likewise transitive on each set
$\mathcal{S_{-\varphi}}$, $\mathcal{S}_{\varphi^-1}$ and
$\mathcal{S}_{-\varphi^{-1}}$, which are the sets of vertices whose
inner product with $v$ is $-\varphi$, $\varphi$$^{-1}$ or
$-\varphi$$^{-1}$, respectively.

This monomial subgroup has two orbits on the set of 30 vertices orthogonal to $v$; these are the six vertices of shape (0, 2, 0, 0) and the 24 of shape (0, 1, $\varphi$, $\varphi$$^{-1}$). 
But if $w$ is the vertex (0, -1, $\varphi$, $\varphi$$^{-1}$) then the reflection $r_w$ (that uses $w$$^{\perp}$ as its hyperplane of reflection) fixes $v$ and fuses the two orbits, so that Stab($v$) is transitive on vertices orthogonal to $v$. 
Similarly, the monomial subgroup $2^3$:$A_3$ has two orbits on the set of 8+12 vertices that have inner product 1 with $v$ (these orbits consist of vertices with shapes (1, 1, 1, 1) and (1, 0, $\varphi$$^{-1}$, $\varphi$)) and once again $r_w$ fuses the orbits. The same is clearly true for the remaining orbit under Stab($v$), which is the set of 20 vertices that have inner product -1 with $v$.

Now, since each of the $\pm60$ vertices of the 600-cell is orthogonal to 15 others, the 600-cell contains 60$\cdot15$/2 = 450 pairs of orthogonal vertices, and these are permuted transitively by 2($A_5$$\times$$A_5$)2 since Stab($v$) is transitive on the $\pm$15 vertices orthogonal to $v$. 
If we choose one of these orthogonal pairs—--say, the vertices ($\pm$2, 0, 0, 0) and (0, $\pm$2, 0, 0)—--then only two of the remaining 58 vertices are orthogonal to both (these are (0, 0, $\pm$2, 0) and (0, 0, 0, $\pm$2), of course) so each of the 450 orthogonal pairs lies in a unique tetrad of mutually orthogonal vertices, i.e., each pair lies in a unique 16-cell. 
Since a 16-cell contains ${4}\choose{2}$ = 6 such pairs, the 600-cell contains exactly 450/6 = 75 16-cells.

But each vertex (pair) in a 24-cell is orthogonal to three others and has inner product $\pm$1 with the remaining eight, so each vertex in a 24-cell lies in exactly one of the three 16-cells in the 24-cell. 
This inner product structure also means that we can complete a 16-cell to a 24-cell only by adding vertices that have inner products $\pm$1 with all four vertices of our given 16-cell. 
For example, if we start with the 16-cell $\{$($\pm$2, 0, 0, 0)$^S$$\}$ then the $\pm$8 other vertices that comprise a 24-cell must have coordinates that are all $\pm$1, and so we can only add the $\pm$8 vertices ($\pm$1, $\pm$1, $\pm$1, $\pm$1) from the 600-cell to obtain a 24-cell.
Thus each 16-cell lies in a unique 24-cell, and a 600-cell therefore contains exactly 75/3 = 25 24-cells. 

Moreover, since an 8-cell is comprised of two 16-cells, where a vertex in either 16-cell has inner product $\pm$1 with each vertex in the other 16-cell, an 8-cell must consist of two of the three 16-cells within the same 24-cell.
Hence there are ${3}\choose{2}$$\cdot$25 = 75 8-cells in a 600-cell. 
And because a 16-cell lies in a unique 24-cell, each 8-cell (being the union of two 16-cells) must lie in a unique 24-cell as well. QED

\medskip

Beyond that, the structure of the 5$\times$5 array described earlier makes it easy to show that \emph{each of the 25 24-cells is disjoint from 8 24-cells and intersects the other 16}. 
Certainly each 24-cell is disjoint from the $4+4$ others that lie in the same row or the same column of the array. 
But given a 24-cell ($i$ $j$), its stabilizer $S_{ij}$ $\cong$ 2($A_4$ $\times$ $A_4$) is transitive on the four rows of the array other than $i$ and the four columns other than $j$, so $S_{ij}$ is transitive on the 16 24-cells whose labels contain neither $i$ nor $j$.
Thus each of these 16 24-cells has the same cardinality of intersection with ($i j$), and this cardinality cannot be 0---otherwise each 24-cell would be disjoint from all the others---so that a 24-cell intersects each of the 16 24-cells outside its row and column of the array.
(In the next section we'll see that each intersection is a regular hexagon.)
But since two 24-cells are disjoint if and only if they lie in the same row or the same column, three or more mutually disjoint 24-cells must all lie in the same exact row or column, so the only partitions of a 600-cell into five disjoint 24-cells are the five rows and the five columns of the array.
Thus Schoute was right all along: there are precisely ten partitions of a 600-cell into five disjoint 24-cells (\textbf{Fact 5}).

And now it's easy to see that elements in the odd half of Aut($H$) interchange the rows and columns of the array (\textbf{Fact 6}). 
For, each row and each column have a single 24-cell in common, whereas any two rows are disjoint (i.e., their $5+5$ 24-cells are distinct) and any two columns are disjoint as well. 
Hence the rows and columns form a bipartite graph, so Aut($H$) permutes Schoute's ten partitions as two pentads: the rows $\{$1, 2, 3, 4, 5$\}$ and the columns $\{$6, 7, 8, 9, X$\}$. 
If $v$ is any vertex of the 600-cell (so that $v$ is labeled (1 $j_1$)(2 $j_2$)(3 $j_3$)(4 $j_4$)(5 $j_5$) by the five 24-cells that contain it) then, as in $\S2$, the reflection $r_v$ fixes each of these five duads, so the reflection either fixes or interchanges the two symbols in each duad.
If $r_v$ mapped rows to rows and columns to columns, then $r_v$ would fix each of the ten symbols (and thus fix each of the $\pm$60 vertices of the 600-cell), so we would have $r_v$ = $\pm$1.
But this can't be true, and the contradiction shows that $r_v$ and all the other elements of Aut($H$) outside 2($A_5$ $\times$ $A_5$) interchange the rows and columns of the array.

\bigskip

\section{The Planes of the 600-cell}

With $\S$7 in mind we describe how pairs of vertices of the 600-cell can interact. 

To begin with, the unit icosian 1$_I$ ($\S$2) is the only vertex that lies in all three 24-cells (16), (27) and (38): the 120 vertices of the 600-cell are conjugate by 2$A_{5R}$ (right multiplication by icosians, which fix $\{$1, 2, 3, 4, 5$\}$ pointwise) and so a vertex labeled (16)(27)(38)(4X)(59) cannot exist, since it would be conjugate to 1$_I$ by only the odd permutation (9X) of $\{$6, 7, 8, 9, X$\}$; and in general, naming any three of a vertex's five duads is sufficient to identify the vertex uniquely. 
Thus the labels of distinct vertices $\{$$\pm$$v$$\}$ and $\{$$\pm$$w$$\}$ cannot share three duads, i.e., pairs of vertices of a 600-cell can lie together in at most two of the 25 24-cells.

If the labels of $\{$$\pm$$v$$\}$ and $\{$$\pm$$w$$\}$ share at least one duad---if $v$ and $w$ are in at least one 24-cell together---then the inner product of $v$ and $w$ is either 0 or $\pm$1. 
If the vertices are orthogonal then they lie in a 16-cell, and therefore in a unique 24-cell, so their labels share only one duad. 
It was shown in $\S$3 that there are 60$\cdot$15/2 = 450 pairs of orthogonal vertices in a 600-cell, with $4 \choose 2$  pairs in each of the 75 16-cells. 
The significance of these 450 pairs can be seen in [10].

If, on the other hand, the labels of $\{$$\pm$$v$$\}$ and $\{$$\pm$$w$$\}$ share two duads (i.e., if they lie in the same two 24-cells) then a third vertex must also lie in both 24-cells in question, since 60/({5$\cdot$4}) = 3 vertex labels contain any two given (disjoint) duads. 
For example, the 24-cells labeled (16) and (27) contain the trio of vertices (16)(27)(38)(49)(5X), (16)(27)(39)(4X)(58) and (16)(27)(3X)(48)(59). 
These $\pm$3 vertices all have inner product $\pm$1 with each other---if any pair were orthogonal, they'd lie in a 16-cell and therefore in a unique 24-cell---so the vertices form a (planar) regular hexagon.

A 24-cell contains 12$\cdot$8/(3$\cdot$2) = 16 hexagons, so \emph{a 24-cell is disjoint from eight 24-cells and intersects each of the other 16 in six vertices that form a hexagon.} 
Since no two of the $\pm$3 vertices in a hexagon have inner product 0, each of the three 16-cells of a 24-cell contains one pair of vertices of a hexagon.
A 600-cell contains 25$\cdot$16/2 = 200 such hexagons, whose significance will be seen in $\S$7.

The final possibility for two vertices $v_1$ and $v_2$ is that their labels have no duads in common---that is, if $v_1$ is labeled (1 $i_1$)(2 $i_2$)(3 $i_3$)(4 $i_4$)(5 $i_5$) and $v_2$  is labeled (1 $j_1$)(2 $j_2$)(3 $j_3$)(4 $j_4$)(5 $j_5$) then each of the $i$'s is a different element of $\{$6, 7, 8, 9, X$\}$ than its $j$ counterpart.
Since 2$A_{5R}$ is transitive on the $\pm$60 vertices of the 600-cell, we can find an icosian $w$ such that $w_R$ fixes $\{$1, 2, 3, 4, 5$\}$ while sending $i_1$ to $j_1$, $i_2$ to $j_2$, etc., so that $w_R$ has order 5 and sends $v_1$ to $v_2$, where we can replace $v_2$ by its negative if necessary.
The orbit of $v_1$ under $w_R$ is therefore five vertices $\{$$v_1$, $v_2$, $v_3$, $v_4$, $v_5$$\}$ that are easily shown to form a (planar) regular pentagon.
Moreover, since $w_R$ fixes $\{$1, ..., 5$\}$ and is transitive on $\{$6, ..., X$\}$, all 25 duads ($i$ $j$) appear among the labels of these five vertices---in other words, \emph{each of the 25 24-cells of the 600-cell contains exactly one vertex of a regular pentagon}.

Naturally the five antipodal vertices $\{$$-v_1$, $-v_2$, $-v_3$, $-v_4$, $-v_5$$\}$ also form a pentagon, and the two pentagons create a regular decagon. 
Since each of the $\pm$60 vertices of a 600-cell has inner product $\varphi$ with 12 others and inner product ${\varphi}^{-1}$ with 12, a 600-cell contains (60$\cdot$24)/(5$\cdot$4) = 72 decagons whose significance will be seen in $\S$7. 
Since each of the 720 edges of the 600-cell is also an edge of a decagon---the two vertices in question have inner product $\varphi$ in either case---\emph{each of the 720 edges of the 600-cell lies on just one of the 72 decagons}.

\bigskip

When we generalize the construction of $\S$2 (which created a
5 $\times$ 5 array from the vertices of the 600-cell), we find that
\emph{the vertices of the 600-cell create an $n \times n$ array for
each of the 600-cell's basic planar shapes}---the square (vertices $v$
and $w$ with $v.w = 0$; $n=5$), the triangle ($v.w$ = $\pm$$1$; $n=10$)
and the pentagon ($v.w$ = $\pm$$\varphi$; $\pm$${\varphi}^{-1}$,
$n=6$.)

We begin with any subgroup $K$ $\leq$ $A_5$ of the icosians. 
(We work projectively for the rest of this section and write $A_5$ instead of 2$A_5$, etc.)
Let $\{$$g_i$$\}$ be a left transversal of $K$ in $A_5$, so that $\{$$g_i$$^{-1}$$\}$ is a right transversal. 
Since our goal is to form  an array of sets $\{$$g^i$$K$$g^{-j}$$\}$ that are distinct for all choices of $i, j$ we require $K$ to be self-normalizing.

The only self-normalizing subgroups of $A_5$ are the normalizers of the Sylow-$p$ subgroups ($p$ = 2, 3 or 5), so we consider the 600-cell from the perspective of these primes. 
To each of the three finite fields $\bf{F}$$_p$ we adjoin the roots of the polynomial $f(x) = x^2-x-1$ in order to incorporate the 600-cell's scalar $\varphi$, which is a root of $f(x)$.
Since $f(x)$ is irreducible mod 2 and mod 3 but splits mod 5, we obtain fields $\bf{F}$$_q$, where $q$ = 4, 9 or 5, depending on $p$.
For the primes 2, 3 and 5 the group $A_5$ contains $q+1$ = 5, 10 and 6 Sylow-$p$ subgroups, respectively. 

Let $N$ denote both the subgroup N($P$) (the normalizer of a Sylow-$p$ subgroup $P$ of $A_5$, where $p$ is 2, 3 or 5) and also the set of vertices we obtain from N($P$) when we consider the $\pm$60 vertices of the 600-cell as the icosians.
Let $\{$$g_i$$\}$ be a left transversal of $N$ in $A_5$, so that $\{$$g_i$$^{-1}$$\}$ is a right transversal.
As in $\S$2, we obtain from $N$ a ($q+1$) $\times$ ($q+1$) array ($q$ is 4, 9 or 5, depending on $p$) where the ($i$, $j$)$^{th}$ entry is the set $g^i$$N$$g^{-j}$ of vertices of the 600-cell. 
Since $N$ is self-normalizing, these sets are all distinct.
Furthermore, since $\{$$g_i$$\}$ is a transversal of $N$ in $A_5$, the first column $\{$$g_i$$N$$\}$ of the array partitions the 120 vertices of the 600-cell into $q+1$ disjoint sets. 
But right multiplication by the icosian $g^{-j}$ is a symmetry of the 600-cell, so each column of the array is a partition; of course the same is also true for each row. 

The case $p$ = 2 was already considered in $\S$3, where $N$ $\cong$ $A_4$ creates a 5 $\times$ 5 array that displays Schoute's ten partitions.
When $p$ = 3, the normalizer $N$ of a Sylow-$p$ subgroup of $A_5$ has order 6; for example, we can take this subgroup (of $A_{5R}$) to be $<$(89X), (67)(89)$>$.
In this case the $\pm$6 vertices $N$ of the 600-cell are the two hexagons (of $\pm$3 vertices)\hfill

\smallskip

{\centering \hspace{0.5cm} $h$ = $\{$(16)(27)(38)(49)(5X), (16)(27)(39)(4X)(58), (16)(27)(3X)(48)(59)$\}$\hspace{0.5cm}  
and $h^{\prime}$ = $\{$(17)(26)(38)(4X)(59), (17)(26)(39)(48)(5X), (17)(26)(3X)(49)(58)$\}$. \hfill}

\smallskip

\noindent Since each vertex of $h$ shares a single duad with each vertex of $h^{\prime}$, the two hexagons are orthogonal, so \emph{the 200 hexagons in a 600-cell form 100 orthogonal pairs that create a 10 $\times$ 10 array.}

When $p$ = 5, the normalizer $N$ of a Sylow-$p$ subgroup of $A_5$ has order 10; for example, we can take this subgroup to be $<$(6789X), (7X)(89)$>$. 
In this case, the $\pm$10 vertices $N$ are the decagons 

$d$ = $\{$(16)(27)(38)(49)(5X), (17)(28)(39)(4X)(56), (18)(29)(3X)(46)(57), \linebreak 
{\centering (19)(2X)(36)(47)(58), (1X)(26)(37)(48)(59)$\}$\linebreak 
and  $d^{\prime}$ = $\{$(16)(2X)(39)(48)(57), (17)(26)(3X)(49)(58), (18)(27)(36)(4X)(59), (19)(28)(37)(46)(5X), (1X)(29)(38)(47)(56)$\}$. \par}

\smallskip

\noindent Since each vertex of $d$ shares a single duad with each vertex of $d^{\prime}$, the two decagons are orthogonal, so \emph{the 72 decagons in a 600-cell form 36 orthogonal pairs that create a 6 $\times$ 6 array.} 

\medskip

Finally, since $L_2(q)$ $\times$ $L_2(q)$ $\cong$ $O_4^+(q)$ ([2], p. xii) and since $L_2(4)$ $\cong$ $L_2(5)$ $\cong$ $A_5$ while $L_2(9)$ $\cong$ $A_6$, each group $O_4^+(q)$.2 ($q$ = 4, 9 or 5) contains a copy of Aut($H$)/$\{$$\pm$1$\}$ $\cong$ ($A_5$ $\times$ $A_5$)2. 
So the arrays can be described in terms of the orthogonal groups $O_4^+(q)$; the case $q=4$ is described in $\S$7 and $q=5$ in [10]. 

\bigskip

\section{The Geometry of $E_8$}
                    
In this section we collect some results about the $E_8$ lattice ([3], p. 120)---or rather, some results about the 8-space $E_8$/2$E_8$ over $\bf{F}$$_2$---that will enable us in $\S$7 to view the $\pm$60 vertices and the 25 24-cells of a 600-cell in a new way. 
We assume some familiarity with the basic results and terminology for vector spaces and quadratic forms over $\bf{F}$$_2$.

From [1, p. 2], the 270 isotropic 4-spaces of $E_8$/2$E_8$ are permuted in two orbits of 135 each by $O_8^{+}(2)$, and given a pair $\{$$V_1$, $V_2$$\}$ of disjoint isotropic 4-spaces (disjoint other than the zero vector, of course) there are 28 others disjoint from both $V_1$ and $V_2$.
These 28 can be labeled by duads on eight letters. 
Since two of the 28 4-spaces are disjoint when their duads share a letter, there are two different ways to complete the pair $\{$$V_1$, $V_2$$\}$ into a maximal set of mutually disjoint 4-spaces (where the next two 4-spaces added towards a maximal set may be labeled by duads (ab) and (ac), and where any two duads that we add must share a letter): (i) we can add the seven 4-spaces $\{$(ab), (ac), (ad),$\ldots$, (ah)$\}$ and complete the original pair $\{$$V_1$, $V_2$$\}$ into a set of nine 4-spaces, or (ii) we can add the three 4-spaces $\{$(ab), (ac), (bc)$\}$ and complete the pair into a set of five 4-spaces.

The automorphism group of the set of nine 4-spaces is the maximal subgroup $A_9$ of $O_8^{+}(2)$ [7, p. 540].
If we relabel the nine 4-spaces $V_1$, $V_2$, (ab), (ac), (ad), $\ldots$, (ah) by 1, 2, $\ldots$, 9, respectively, then the remaining 126 4-spaces in the 135-orbit can be labeled by the ${9}\choose{5}$ = 126 pentads of the nine symbols, according to which five of the nine 4-spaces they intersect: any pair of the 135 4-spaces are either disjoint or intersect in a 2-space (which contains three non-zero vectors), and since our nine disjoint 4-spaces contain all 135 = 9$\cdot$15 isotropic vectors among them, each of the 126 4-spaces has its 15 non-zero vectors allocated among exactly five of the nine.
But then the five disjoint 4-spaces of case (ii) above are $\{$1, 2, 3, 4, (56789)$\}$, and the subgroup of $A_9$ that stabilizes this pentad of 4-spaces is ($A_4$$\times$$A_5$)2.

The only 4-spaces that intersect all four of the 4-spaces $\{$1, 2, 3, 4$\}$ are $\{$(12345), (12346), (12347), (12348), (12349)$\}$. 
But these five 4-spaces are also mutually disjoint$\footnote{Two 4-spaces are disjoint when their labels (by one of 9 monads or one of 126 pentads) intersect with even cardinality.}$, and each one intersects all five 4-spaces $\{$1, 2, 3, 4, (56789)$\}$ in a 2-space of three (non-zero) vectors; in other words, the 15 vectors in any of the ten 4-spaces consist of three vectors from each 4-space in the other pentad, and so both pentads contain the same 5$\cdot$15 = 75 isotropic vectors. 
Since $O_8^{+}(2)$ is transitive on tetrads of disjoint 4-spaces there must exist a group element that interchanges the two pentads, so the automorphism group of the two pentads is the maximal subgroup ($A_5$$\times$$A_5$)$2^2$ of $O_8^{+}(2)$ [2, p. 85].

Furthermore, the intersection pattern of the two pentads creates the same 5$\times$5 array as the 25 24-cells of the 600-cell. 
We can arrange the 75 isotropic vectors common to both pentads into a 5$\times$5 array of 25 2-spaces (each having three non-zero vectors), where the five rows are headed by the first pentad $\{$$V_i$$\}$, the five columns are headed by the second pentad $\{$$W_j$$\}$, and the $(ij)^{th}$ entry is the 2-space $\ell_{ij}$ that contains the three vectors that lie in both $V_i$ and $W_j$, as in Figure 1.

\bigskip

\hspace{3cm{}}$W_{1}$\hspace{0.95cm}$W_{2}$\hspace{0.95cm}.\hspace{1cm}.\hspace{1cm}$W
_{5}$\hfill

\medskip

\hspace{1.7cm}$V_1$\hspace{0.95cm}$\ell_{11}$\hspace{1cm}$\ell_{12}$\hspace{1cm}.\hspace{1cm}.\hspace{1cm}$\ell_{15}$\hfill

\hspace{3.1cm}.\hspace{1.5cm}.\hspace{3.5cm}.\hfill

\hspace{1.7cm{}}$V_5$\hspace{0.9cm}$\ell_{51}$ \hspace{1cm}$\ell_{52}$\hspace{1cm}.\hspace{1cm}.\hspace{1cm}$\ell_{55}$\hfill

\begin{center}\textit{Figure 1. Two pentads of isotropic 4-spaces and their 25 lines (2-spaces) of intersection}\end{center}

The identification of the 25 isotropic 2-spaces in Figure 1 with the 25 24-cells of a 600-cell will be made apparent in $\S$7.

\bigskip

\section{Changing the Norms of Vectors}

In this section we define maps from quadratic fields to $\mathbb{Q}$ in a variety of ways, which leads a variety of embedded $H_4$ polytopes. Our main objective will be to embed a 600-cell in the $E_8$ lattice.

\medskip

The coordinates and inner products for $H_4$ polytopes lie in the golden field $\mathbb{Q}$($\sqrt{5}$), making it impossible to (legitimately) embed $H_4$ polytopes in the $E_8$ lattice, whose inner products lie in $\mathbb{Z}$. 
But it was shown in [15] that a “reduced norm” for the 120 vertices of a 600-cell $H$ identifies them as 120 of the 240 root vectors of $E_8$, and the remaining 120 root vectors come from the same norm-reducing formula applied to the scaled 600-cell $\varphi$$H$.
The reduction formula is obtained by composing the standard inner product in $\mathbb{Q}$($\sqrt{5}$) with a modified trace map from $\mathbb{Q}$($\sqrt{5}$) to $\mathbb{Q}$, and the trace map of [15] identifies a vertex $v$ in $H$ or $\varphi$$H$ of norm \textit{a} + \textit{b}$\sqrt{5}$ (\textit{a, b} $\in$ $\mathbb{Q}$) with a vector v $\in$ $E_8$ of reduced norm \textit{a} - \textit{b}. ([15], p. 57)
 
But [16] uses a different trace map to embed the 600-cell in $E_8$
([16], p. 158), which raises the question of how many different norm
reductions are possible. 
The answer is as follows: \hfill\break

$\textbf{Proposition}$: Let $\mathbb{Q}$($\sqrt{n}$) ($n$ $>$ 0) be a real quadratic field. 
Then for $m$ $\in$ $\mathbb{Q}$ there exists a norm reduction from $\mathbb{Q}$($\sqrt{n}$) to $\mathbb{Q}$ that maps $\sqrt{n}$ to $m$, i.e., 
\begin{center}\textit{a} + \textit{b}$\sqrt{n}$ $\rightarrow$  \textit{a} + \textit{bm}  $\quad$   (\textit{a, b} $\in$ $\mathbb{Q}$)\end{center}
if and only if $|$$m$$|$ $<$ $\sqrt{n}$.\hfill\break

\smallskip Proof: The norm of a typical scalar \textit{x} + $y$$\sqrt{n}$ in $\mathbb{Q}$($\sqrt{n}$) (where \textit{x, y} $\in$ $\mathbb{Q}$) must reduce to a positive number in $\mathbb{Q}$. 
But under the map that reduces $\sqrt{n}$ to $m$, this scalar's norm of

\hspace{1.7cm}($x$ + $y$$\sqrt{n}$$)^2$ = $x^2$ + 2$xy$$\sqrt{n}$ + $n$$y^2$ \hspace{0.2cm} reduces to 

\hspace{4cm}$x^2$ + 2$xym$\hspace{0.25cm} + $n$$y^2$ = ($x$+$my$)$^2$ + ($n-m^2$)$y^2$.

When ($x$, $y$) = ($-m$, 1) this norm is $n-m^2$, so we must have $|$$m$$|$ $<$ $\sqrt{n}$. 
But in that case---whenever $|$$m$$|$ $<$ $\sqrt{n}$---if a vector in a space over $\mathbb{Q}$($\sqrt{n}$) has coordinates $\{$($x_i$+$y_i$$\sqrt{n}$)$\}$, we can split each coordinate  into two coordinates ($x_i$+$my_i$, $\sqrt{n-m^2}$$y_i$) to obtain a norm reduction that maps $\sqrt{n}$ $\rightarrow$ $m$. QED

\medskip

In the particular case of lattices over $\mathbb{Q}$($\sqrt{5}$) we obtain $\mathbb{Q}$-lattices under the norm-reduction map 

\begin{center}\textit{a} + \textit{b}$\sqrt{5}$ $\rightarrow$  \textit{a} + \textit{bm}  $\quad$   ($m$ $\in$ $\mathbb{Q}$)\end{center}
for any $|$$m$$|$ $<$ $\sqrt{5}$, and for an integral lattice over $\mathbb{Z}$[$\varphi$] we obtain an integral $\mathbb{Z}$-lattice for $m$ $\in$ $\{$0, 
$\pm$1, $\pm$2$\}$, where the coordinates for these values of $m$ are given by ($x_i$+$y_i$$\sqrt{5}$)
$\rightarrow$ $\{$($x_i$, $y_i$$\sqrt{5}$), ($x_i$$\pm$$y_i$, 2$y_i$), or ($x_i$$\pm$2$y_i$, $y_i$)$\}$
respectively. 
Note that the natural norm may need to be doubled when $m$ $\in$ $\{$0, $\pm$2$\}$ since $\varphi$ = $\frac{1}{2}$(1+$\sqrt{5}$) is an integer in the golden field that reduces to $\frac{1}{2}$, $\frac{3}{2}$, $-\frac{1}{2}$ for $m$ = 0, 2, $-2$, respectively.

\medskip

The following are a few examples of how the Proposition may be used to transform $H_4$ polytopes into lattices over $\mathbb{Q}$. 

Example 1: We begin with the 600-cell. 
Then the norm-changing maps
\begin{center}$a+b\sqrt{5}$ $\rightarrow$ $a-b$, $\quad$ $a+b\sqrt{5}$ $\rightarrow$ $a+b$ $\quad$ ($m = -1, 1$, resp.)\end{center}
both embed the 600-cell as 120 of the 240 norm 2 vectors in $E_8$. 
Note that [15] used $m = -1$ (and mapped $\varphi$ = $\frac{1}{2}$$(1+\sqrt{5}$) to 0) whereas [16] used $m$ = 1 
(and mapped $\varphi$$^{-1}$ = $\frac{1}{2}$(-1+$\sqrt{5}$) to 0.)

\medskip

Example 2: Using the 600-cell again, we can choose $m$ = 0 instead, so that $a+b\sqrt{5}$ $\rightarrow$ $a$. 
In this case, if we double the natural inner product then we obtain a rank 8 $\mathbb{Z}$-lattice L generated by $\pm60$ norm 4 vectors, each of which has inner product\hfill

\hspace{2.1cm}$\pm$4 with 1 pair of vectors;\hfill

\hspace{2cm}	$\pm$2 \hspace{0.7cm}20;\hfill

\hspace{2cm}    $\pm$1   \hspace{0.7cm}24; and \hfill

\hspace{2.1cm} 0  \hspace{0.9cm}15.\hfill\break

Now, we can find four vertices $\{$$v_1$, $v_2$, $v_3$, $v_4$$\}$ of the 600-cell $H$ so that all $\pm$60 vertices of $H$ are $\mathbb{Z}$[$\varphi$]-linear combinations of the $v_i$’s. 
It’s easily verified that the $\mathbb{Z}$[$\varphi$]-lattice H
generated by these four vertices is self-dual, so if we let 
$\{$$w_1$, $w_2$, $w_3$, $w_4$$\}$ denote the dual basis of H 
(i.e., inner products satisfy $v_i$$\cdot$$w_j$ = $\delta$$_{ij}$) 
then the $w_j$’s are $\mathbb{Z}$[$\varphi$]-linear combinations of the $v_i$’s.
When we apply our norm reduction procedure to the lattice H (still using $m$ = 0) each coordinate $\{$$a + b$$\sqrt{5}$$\}$ of a vector $v$ in H splits into two coordinates $\{$$a, b$$\sqrt{5}$$\}$ of the norm-reduced vector $v_L$ of the $\mathbb{Z}$-lattice L, as above.

But L is generated by the eight vectors $\{$($v_i)$$_L$,
($\varphi$$v_i)$$_L$$\}$, and since we have
($v_i$)$_L$$\cdot$($w_j$)$_L$ = 2$\delta_{ij}$ (recall that the natural
inner product has been doubled), and since we also have
($\varphi$$v_i$)$_L$$\cdot$($w_j$)$_L$ = $\delta_{ij}$ = 
($v_i$)$_L$$\cdot$($\varphi$$w_j$)$_L$ and
($\varphi$$v_i$)$_L$$\cdot$($\varphi$$w_j$)$_L$ = 3$\delta_{ij}$, the
basis dual to $\{$($v_i)$$_L$, ($\varphi$$v_i)$$_L$$\}$ is given by
$\{$$\frac{1}{5}$($3w_i$ – $\varphi$$w_i$)$_L$, $\frac{1}{5}$($-w_i$ +
2$\varphi$$w_i$)$_L$$\}$. 

Thus d(L) = $5^4$, and L must be the unique 8-dimensional rootless even lattice $\Lambda$ with d($\Lambda$) = $5^4$ [8, p. 3].

\medskip

Example 3: The 720 vertices of the rectified 600-cell (whose vertices are the midpoints of each of the 720 edges of the 600-cell) can be given coordinates of shapes 
(0, 0, 2$\varphi$, 2$\varphi$$^2$), (1, 1, $\varphi^3$, $\varphi^3$), (0, 1, $\varphi$, 1+3$\varphi$), (0, $\varphi^2$, $\varphi^3$, 2+$\varphi$), (1, $\varphi$, 2$\varphi^2$, $\varphi^2$) and ($\varphi$, $\varphi^2$, 2$\varphi$, $\varphi^3$).
These vertices have norm $20 + 8\sqrt{5}$, and by choosing $m=-2$ (which maps $a+b\sqrt{5}$
$\rightarrow$ $a - 2b$) the rectified 600-cell is embedded as 720 norm 4 vectors of $E_8$.
In fact, it's easy to show that the 2160 norm 4 vectors of $E_8$ consist of five embedded $H_4$ polytopes: two 600-cells (of 120 vertices each); two 120-cells (of 600 vertices each); and a rectified 600-cell of 720 vertices.

\section{Aut($H_4$) $\rightarrow$ $O_4^{+}$(4).2}

The variety of norm-reducing maps (and their resulting embeddings) means that a lattice---or polytope---can be viewed in a number of new ways. 
This is the case when we embed 600-cells $H$ and $\varphi$$H$ in $E_8$ and project them to $E_8$/$2E_8$: the 8-space over $\bf{F}$$_2$ employs the geometry of the 600-cell to become a 4-space over $\bf{F}$$_4$ (\textbf{Fact 10}), as we now demonstrate. 

Let $\varepsilon$: $H$ $\rightarrow$ $E_8$ denote the embedding of $H$ in $E_8$ used in [15] (i.e., the norm reduction map of $\S$6, Example 1 that uses $m$ = $-1$), and similarly let $\varepsilon$$_{\varphi}$: $\varphi$$H$ $\rightarrow$ $E_8$ denote the embedding of $\varphi$$H$. 
We now construct a map $\Phi$: $E_8$ $\rightarrow$ $E_8$ that sends the each of the 120 root vectors embedded from $H$ onto its counterpart---the corresponding root vector embedded from $\varphi$$H$---by defining $\Phi$($\varepsilon$($v$)) = $\varepsilon$$_{\varphi}$($\varphi$$v$) for any vertex $v$ $\in$ $H$.
Note that since Aut($H$) $\cong$ 2($A_5$$\times$$A_5$)2 permutes $H$ and $\varphi$$H$ exactly the same way, $\Phi$ commutes with all of 2($A_5$$\times$$A_5$)2 $\leq$ Aut($E_8$).

In terms of coordinates, if $\alpha$ = \textit{a} + \textit{b}$\sqrt{5}$ (\textit{a, b} $\in$ $\mathbb{Q}$) is one of the four coordinates of a vertex $v$ $\in$ $H$ then $\alpha$ splits into two coordinates ($a-b$, $2b$) when $H$ is embedded in $E_8$, as in $\S$6. 
Likewise, the scaled coordinate $\varphi$$\alpha$ = $\frac{1}{2}$[($a + 5b) + (a + b$)$\sqrt{5}$)] (the corresponding coordinate of $\varphi$$v$) splits into two coordinates ($2b, a + b$) under the embedding $\varepsilon$$_{\varphi}$. 
Setting $x$ = $a-b$ and $y = 2b$, we find that $\Phi$ maps a split coordinate pair ($x$, $y$) to ($y$, $x+y$).

Since the vertices $v_0$ = (2, 0, 0, 0) and $\varphi$$v_0$ = (1+$\sqrt{5}$, 0, 0, 0) of $H$ embed in $E_8$ as the orthogonal root vectors (2, 0, $\dots$, 0) and (0, 2, 0, $\dots$, 0) and since $\Phi$ commutes with all of 2($A_5$$\times$$A_5$)2 $\leq$ Aut($E_8$), every root vector $\varepsilon$($v$) embedded from $H$ is orthogonal to its counterpart $\varepsilon$$_{\varphi}$($\varphi$$v$).
Thus v+$\Phi$(v) is a norm 4 vector of $E_8$ whenever v $\in$ Im($\varepsilon$).

In addition, since the 120 root vectors embedded from $H$ generate $E_8$, $\Phi$ extends to a (well-defined) endomorphism of $E_8$, and so we now have $\Phi$(g(v)) = g($\Phi$(v)) for all g $\in$ 2($A_5$ $\times$ $A_5$)2 and all vectors v $\in$ $E_8$. 
Since $\Phi$ maps a pair ($x$, $y$) to ($y$, $x+y$), $\Phi$$^2$ maps ($x$, $y$) to ($x+y$, $x+2y$); thus $\Phi$ retains the defining property of $\varphi$ = $\frac{1}{2}$(1+$\sqrt{5}$):

\medskip

\hspace{1.5cm} (i): $\Phi$$^2$ = $\Phi$ + 1, \hspace{0.5cm}so that \hspace{0.5cm} (ii): $\Phi$$^3$ = 2$\Phi$ + 1

\medskip

If $\overline{\Phi}$ denotes the endomorphism of $E_8$/2$E_8$ induced by $\Phi$, then by (ii), $\overline{\Phi}$$^3$ = 1, and by (i) the orbit of any vector $v$ $\in$ $E_8$/2$E_8$ under $\overline{\Phi}$ is the set $\{$$v$, $\overline{\Phi}$($v$), $v$+$\overline{\Phi}$($v$)$\}$--that is, the three non-zero vectors in a 2-space over $\bf{F}$$_2$.
But then $\overline{\Phi}$ is an automorphism that transforms our 8-space over $\bf{F}$$_2$ into a 4-space over $\bf{F}$$_4$: we have an element $\omega$ (=  $\overline{\Phi}$) where $\omega$$^3$ = 1, and the orbit of any vector $v$ under $\omega$ is the set $\{$$v$, $\omega$$v$, $\overline{\omega}$$v$ (=$v$+$\omega$$v$)$\}$.
In addition, as we'll see below, the quadratic form on the 8-space over $\bf{F}$$_2$ gives rise to a quadratic form on our 4-space over $\bf{F}$$_4$, and since $\Phi$ commutes with 2($A_5$ $\times$ $A_5$)2 we also have g($\omega$$v$) = $\omega$g($v$) for all g $\in$ ($A_5$ $\times$ $A_5$)2 $\cong$ $O_4^{+}$(4).2.

If $v$ is the image in $E_8$/2$E_8$ of a root vector embedded from $H$, then the three (non-zero) vectors $\{v$, $\overline{\Phi}(v)$, $v$+$\overline{\Phi}$($v$)$\}$ of a 2-space consist of the non-isotropic vector $v$ that came from $H$, its (non-isotropic) counterpart $\overline{\Phi}(v)$ embedded from $\varphi$$H$, and the isotropic vector $v+\overline{\Phi}$($v$) (recall that $\varepsilon$($v_0$) + $\varepsilon$$_\varphi$($\varphi$$v_0$) is a norm 4 vector of $E_8$ for any vertex $v_0$ of $H$, as above.) 
Thus 180 of the 255 non-zero vectors of $E_8$/2$E_8$ form 60 2-spaces that can be labeled by the $\pm$60 vertices $\{\pm v\}$ of $H$, and these 2-spaces account for all 60 + 60 non-isotropic vectors and 60 of the 135 isotropic vectors of $E_8$/2$E_8$.

The 75 remaining isotropic vectors form the 25 2-spaces of Figure 1 ($\S$5), where each 2-space $\ell_{ij}$ has a distinct stabilizer $S_{ij}$ of index 25 in ($A_5$$\times$$A_5$)2 $\leq$ Aut($E_8$/2$E_8$).
Since $\overline{\Phi}$ commutes with all of ($A_5$$\times$$A_5$)2, we have $S_{ij}$($\overline{\Phi}$($\ell_{ij}$)) = $\overline{\Phi}$($S_{ij}$($\ell_{ij}$)) = $\overline{\Phi}$($\ell_{ij}$), so that $\overline{\Phi}$ stabilizes each 2-space and therefore permutes its three (non-zero) vectors.

But the isomorphism between subgroups ($A_5$$\times$$A_5$)2 of both Aut($E_8$/2$E_8$) and Aut($H$)/$\{$$\pm1$$\}$ means that each stabilizer $S_{ij}$ $\leq$ Aut($E_8$/2$E_8$) corresponds to a specific subgroup of index 25 in Aut($H$)/$\{$$\pm1$$\}$.
Since these subgroups lift in Aut($H$) to the stabilizers of the 25 24-cells of the 600-cell, each 2-space of Figure 1 corresponds to a specific 24-cell of $H$.

So we have obtained \textbf{Fact 10}: embedding the 600-cell $H$ in $E_8$ turns the 8-space over $\bf{F}$$_2$ into a 4-space over $\bf{F}$$_4$ whose 85 (projective) points are labeled by the 60 vertex pairs and the 25 24-cells of $H$. 

\medskip
 
Our 4-space over $\bf{F}$$_4$ derives a quadratic form Q$_{\omega}$ from the quadratic form Q on the 8-space over $\bf{F}$$_2$. 
If $v$ is any vector of $E_8$/2$E_8$ then we define Q$_{\omega}$($v$) = Q($v$) + $\omega$Q($\omega$$v$) + $\overline{\omega}$Q($\overline{\omega}$$v$); a simple calculation shows that Q$_{\omega}$ is a quadratic form.
If the vectors $\{$$v$, $\omega$$v$, $\overline{\omega}$$v$$\}$ are all isotropic then Q$_{\omega}$($v$) = 0 (so the 25 points labeled by 24-cells are singular); if $v$ is isotropic but $\omega$$v$ and $\overline{\omega}v$ are non-isotropic then Q$_{\omega}$($v$) = 1; and if $v$ is non-isotropic then Q$_{\omega}$($v$) = $\overline{\omega}$ if $v$ is the image in $E_8$/2$E_8$ of a vector embedded from $H$ and Q$_{\omega}$($v$) = $\omega$ when $v$ came from $\varphi$$H$. 

Naturally the $\bf{F}$$_2$ quadratic form Q can be recovered from Q$_{\omega}$ using a trace formula: we have Q($v$) =

\begin{center}$Tr_{\bf{F_4/F_2}}$(Q$_{\omega}$($v$)) = $\left\{
  \begin{array}{rcr}
      0, & $\textrm{if Q$_{\omega}$($v$) $\in$ $\{$0, 1$\}$}$ &\\
      1, & $\textrm{if Q$_{\omega}$($v$) $\in$ $\{$$\omega$, $\overline{\omega}$$\}$}$  \\
  \end{array}
\right.$
\end{center}

\bigskip

The isomorphism Aut($H_4$)/$\{$$\pm$1$\}$ $\cong$ $O_4^{+}$(4).2 means that certain types of geometric objects are favored in the 600-cell; these correspond to the $\frac{4^{4}-1}{4-1}$ = 85 1-spaces ("points")---and by duality the 85 3-spaces ("planes")---plus the $\frac{85\cdot84}{5\cdot4}$ = 357 2-spaces ("lines") in a 4-space over $\bf{F}$$_4$. 
Each line contains $\frac{4^{2}-1}{4-1}$ = 5 points and each plane contains $\frac{4^{3}-1}{4-1}$ = 21 points.

\bigskip

The following proposition will help determine the structure of the lines and planes in terms of their constituent points.\hfill\break

$\textbf{Proposition}$: The subgroup $S_C$ of Aut($H$) that stabilizes a 24-cell $C$ in a 600-cell $H$ is isomorphic to 2($A_4$$\times$$A_4$)2, and the subgroup $S_v$ that fixes a vertex $v$ of $H$ is isomorphic to 2$\times$$A_5$ (\textbf{Fact 4}).

\smallskip Proof: Aut($H$) $\cong$ 2($A_5$$\times$$A_5$)2 is transitive on the 25 24-cells in their 5$\times$5 array, so the first statement is obvious.

Since Aut($H$) is transitive on the 120 vertices of the 600-cell, $S_v$ has order 14,400$/$120 = 120. Furthermore, since $v$ is negated by $r_v$ (the reflection that uses $v^{\perp}$ as its hyperplane), $-r_v$ is an element of order 2 that fixes $v$. But $-r_v$ is an odd permutation of $\{$1, 2, ..., 9, X$\}$ that lies in the center of $S_v$, so that $S_v$ $\cong$ $<$$-r_v$$>$$\times$$S_v^e$, where $S_v^e$ is the even subgroup of index 2 in $S_v$. 

Recall ($\S$2) that the vertex $v$ is labeled by five duads (1 $j_1$)(2 $j_2$)(3 $j_3$)(4 $j_4$)(5 $j_5$), and that left multiplications $A_{5L}$ achieve all even permutations of $\{$1, 2, 3, 4, 5$\}$ while right multiplications $A_{5R}$ achieve all even permutations of $\{$6, 7, 8, 9, X$\}$. 
Thus every even permutation of the five duads (that is, every even permutation of the five 24-cells that contain the vertex $v$) can be achieved by an element $w_L$$w^{\prime}$$_R$ of Aut($H$), where $w$ is chosen so that $w_L$ permutes $\{$1, 2, 3, 4, 5$\}$ exactly as the duads are to be permuted, and $w^{\prime}$ is chosen so that $w^{\prime}$$_R$ sends each of $\{$6, 7, 8, 9, X$\}$ back to its original duad partner. 
Since $w_L$$w^{\prime}$$_R$ preserves the set of five duads, $v$ is either fixed or negated. 
If $v$ is negated we replace $w_L$$w^{\prime}$$_R$ by $-w_L$$w^{\prime}$$_R$, so that $w_L$$w^{\prime}$$_R$ $\in$ $S_v^e$ and thus $S_v^e$ maps onto $A_5$.
Since the order of $<$$-r_v$$>$$\times$$S_v^e$ is 120, we must have $S_v^e$ $\cong$ $A_5$ and thus $S_v$ $\cong$ 2$\times$$A_5$. QED

\bigskip

The identification of the 85 points as 25 24-cells ($i j$) and $\pm$60 vertices  (1 $j_1$)(2 $j_2$)(3 $j_3$)(4 $j_4$)(5 $j_5$) makes it easy to determine the lines and planes of our 4-space, simply by observing how $S_C$ and $S_v$ permute these 85 objects via their labels.
For example, the stabilizer $S_C$ of a 24-cell ($i j$) permutes the other 24-cells in orbits of sizes 8 and 16, where the orbit of size 8 consists of the 24-cells whose labels contain either $i$ or $j$. 
Furthermore, this subgroup permutes the 60 vertices in orbits of sizes 12 and 48, according to whether or not ($i j$) is one of the five duads of the vertex in question---that is, whether or not the vertex is contained in the 24-cell $C$.

The vertex stabilizer $S_v$ permutes the 25 24-cells in orbits of sizes 5 and 20, where the orbit of size 5 consists of the five duads that label $v$ (i.e., the five 24-cells that contain $v$.) And $S_v$ permutes the 60 vertices in orbits of sizes 1 + 15 + 20 + 12 + 12, depending on whether the inner product with $v$ (modulo $\{$$\pm$1$\}$) is 2, 0, 1, $\varphi$ or $\varphi$$^{-1}$, respectively, as in the proof of the Theorem in $\S$2. 
(Equivalently, the orbits are determined by the conjugacy class of the $A_5$ permutation of $\{$6, 7, 8, 9, X$\}$ that takes the five duads of $v$ to the five duads of the other vertex.)

It's now a simple exercise to determine the different types of lines and planes in our 4-space over $\bf{F}$$_4$. 
For example, let $\ell$ be a line (of five points) generated by a vertex-labeled point and a 24-cell-labeled point, where the vertex $v$ is contained in the 24-cell $C$. 
The stabilizer $S_v$ of the vertex $v$ is isomorphic to 2$\times$$A_5$, so the stabilizer of both $v$ and the 24-cell $C$ is a subgroup $S_{C, v}$ = $S_C$ $\cap$ $S_v$ $\cong$ 2$\times$$A_4$. 
As above, $S_C$ permutes the 25 24-cells in orbits of sizes 1 + 8 + 16, and $S_v$ permutes the 25 in orbits of sizes 5 and 20.

The 5-orbit of $S_v$ decomposes under $S_{C, v}$ as 1 + 4 since $S_{C, v}$ fixes $C$ and is transitive on the four other 24-cells that contain $v$. 
The 20-orbit decomposes as 8 + 12 since $S_{C, v}$ is transitive on the 8 of 20 duads (24-cells) that share one symbol with $C$, and transitive on the 12 that contain neither symbol.
In particular, the line $\ell$ cannot contain another point labeled by a 24-cell, since $S_{C, v}$ already fixes two of the line's five points (the ones labeled by $v$ and by $C$) and so must permute the other three; but $S_{C, v}$ permutes the 25 24-cells as 1 + 4 + 8 + 12, so (other than the fixed point of $C$ itself) the orbit sizes are too large.

So the other three points on the line must all be labeled by vertices, and we now consider the orbits of $S_{C, v}$ on the 60 vertices. 
The orbit of size 15 under $S_v$ consists of the 15 vertices orthogonal to $v$, three of which lie in the 24-cell $C$ and 12 of which lie in the other four 24-cells that contain $v$. 
Since $S_{C, v}$ $\cong$ 2 $\times$ $A_4$ is transitive on these four 24-cells, the orbit of size 15 splits as 3 + 12.

The orbit of size 20 under $S_v$ splits into orbits of sizes 8 + 12 under $S_C$: the 20 vertices in this orbit have inner product 1 with $v$, and 8 of these 20 lie in $C$ (and are permuted transitively by $S_{C, v}$) while the other 12 are split among  the other four 24-cells that contain $v$.
These 12 are also permuted transitively. 

Finally, the orbits of size 12 under $S_v$ must remain intact under $S_{C, v}$, since the stabilizer in $S_v$ of one of the 12 vertices has order 120/12 = 10; but the order of $S_{C, v}$ isn't divisible by 5, so its stabilizer subgroup can have order at most 2, resulting in the same orbit size of 24/2 = 12.

So only one orbit of $S_{C, v}$ contains as few as three points, and these points are labeled by the three vertices of $C$ that are orthogonal to $v$. 
Thus the line $\ell$ consists of a single point labeled by a 24-cell $C$ and four points labeled by the vertices of a 16-cell in $C$, and so there are 75 = 25$\cdot$3 such lines in our geometry.

Similar arguments determine the other three types of lines.
The two types of planes are especially simple to identify, since by duality their stabilizers are the 25 + 60 subgroups $S_C$ and $S_v$.

\medskip

Thus we find our geometry has:

\medskip

\noindent(A) 85 = 60 + 25 points that correspond to:

(i)	The $\pm$60 vertices of a 600-cell

(ii) The 25 24-cells of a 600-cell

\medskip

\noindent(B)	85 = 60 + 25 planes of 21 points each, corresponding to:

(i)$^{\prime}$  The $\pm$60 vertices. The 21 = 1 + 15 + 5 points correspond to a vertex, the 15 vertices orthogonal to it, and the 5 24-cells that contain it.

(ii)$^{\prime}$   The 25 24-cells. The 21 = 1 + 8 + 12 points correspond to a 24-cell, the 8 24-cells disjoint from it, and the 12 vertices it contains.

\medskip

Finally, the 357 lines honor Schoute's ten partitions alongside the three planar shapes ($\S$4) that a 600-cell's vertices create:

\medskip

\noindent(C)	357 = 10 + 72 + 75 + 200 lines of 5 points each, corresponding to:

(iii)	Schoute’s ten partitions of the 120 vertices of a 600-cell. 
The 0+5 points correspond to a set of five disjoint 24-cells. (Schoute's ten lines are totally singular.)

(iv)	$\pm$72 pentagons of 5 vertices each. 
The 5+0 points correspond to the 5 vertices of a pentagon.   ($<$v.w$>$ = $\varphi$ or $\varphi^{-1}$, 72 = $\frac{60\cdot24}{5\cdot4}$).

(v)	75 16-cells of 4 mutually orthogonal vertices each. 
The 4+1 points correspond to the four vertices of a 16-cell, and the 24-cell that contains it. ($<$v.w$>$ = 0, 75 = $\frac{60\cdot15}{4\cdot3}$).

(vi)	$\pm$200 triangles of 3 vertices each. 
The 3+2 points correspond to the triangle of vertices whose labels contain both duads ($i$ $j$) and ($i'$ $j'$), plus the two related 24-cells ($i$ $j'$) and ($i'$ $j$). ($<$v.w$>$ = 1, 200 = $\frac{60\cdot20}{3\cdot2}$).\\

$\vspace{1 cm}$

\begin{center}$\bf{References}$\end{center}

\medskip

[1] Butler, M. et al, \textit{The Unknown Subgroup of Aut($E_8$)}, ArXiv 1709.05532, 2017\hfill

[2] Conway, J.H.C. et al., \textbf{ATLAS of Finite Groups}, Clarendon Press, Oxford, 1985

[3] Conway, J.H.C. $\&$ and Sloane, N.J.A, \textbf{Sphere Packings, Lattices and Groups}, $3^{rd}$ edition, Springer-Verlag, 1999

[4] Coxeter, H.S.M., $\textit{Regular Compound Polytopes in More than Four Dimensions}$, \textbf{J. Math Phys. 12} (1933), 334-345\hfill

[5] Coxeter, H.S.M., $\textbf{Regular Polytopes}$, $3^{rd}$ edition, Dover, 1973\hfill

[6] Du Val, P., $\textbf{Homographies, Quaternions and Rotations}$, Clarendon Press, Oxford, 1964\hfill

[7] Dye, R.H., \textit{The Simple Group FH(8, 2) of Order $2^{12}3^55^27$ and the Associated Geometry of Triality}, \textbf{Proc LMS 18} (1968), 521-562\hfill

[8] Griess, R.L. $\&$ Lam, C.H., \textit{A Moonshine Path for 5A and Associated Lattices of Ranks 8 and 16}, \textbf{Journal of Algebra 331} (2011), 338-361; ArXiv 1006.390\hfill

[9] Jacobson, N., \textbf{Basic Algebra I}, $2^{nd}$ edition, Dover, 2005\hfill

[10] Johnson, D., \textit{$H_4$ over $\bf{F}$$_5$}, ArXiv\hfill

[11] Littlewood, D.E., \textit{The Groups of the Regular Solids in $n$ Dimensions}, \textbf{Proc LMS s2-32} (1931), 10-20\hfill

[12] Schl\"afli, L., \textit{Réduction d’une intégrale multiple, qui comprend l’arc de cercle et l’aire du triangle sphérique comme cas particuliers}, \textbf{Journal de Mathématiques Pures et Appliquées (1), 20} (1855), 359-394

[13] Schoute, P.H., \textbf{Mehrdimensionale Geometrie}, vol. 2, Leipzig 1905\hfill

[14] Somerville, D.M.Y., \textbf{Introduction to the Geometry of n Dimensions}, Methuen $\&$ Co., London, 1929\hfill

[15] Tits, J., \textit{Quaternions over $\mathbb{Q}$($\sqrt{5}$), Leech's Lattice and the Sporadic Group of Hall-Janko}, \textbf{Journal of Algebra 63} (1980), 56-75\hfill

[16] Wilson, R.A., \textit{The Geometry of the Hall-Janko Group as a Quaternionic Reflection Group}, \textbf{Geometriae Dedicata 20} (1986), 157–173

\bigskip

Contacts: {tomme.denney32@gmail.com, babyshay12@iclould.com, jdejaneke@yahoo.com,  tiannarobinson71401@gmail.com,
majidhbutler3509@gmail.com, nisha.claiborne@gmail.com\\}

\end{document}